\documentclass[12pt,reqno]{article}

 \usepackage{amsmath, amsthm, amsfonts}
 \usepackage{mathrsfs, amssymb}
 \usepackage{graphicx, color, bm, cases}

 \numberwithin{equation}{section}

\newtheorem{thm}{Theorem}[section]

\newtheorem{lemma}[thm]{Lemma}

\newtheorem{prop}[thm]{Proposition}
\newtheorem{rem}[thm]{Remark}

\newtheorem{con}[thm]{Condition}
\renewcommand {\theequation}{\thesection.\arabic{equation}}

\def\ga{{\gamma}}

\def\<{\left<}\def\>{\right>}\def\({\left(}\def\){\right)}

\newfam\msbmfam\font\tenmsbm=msbm10\textfont
\msbmfam=\tenmsbm\font\sevenmsbm=msbm7

\scriptfont\msbmfam=\sevenmsbm

 \def\beqlb{\begin{eqnarray}}\def\eeqlb{\end{eqnarray}}
 \def\beqnn{\begin{eqnarray*}}\def\eeqnn{\end{eqnarray*}}
 \def\d{{\mbox{\rm d}}}\def\e{{\mbox{\rm e}}}

 \def\mcr{\mathscr}\def\mbb{\mathbb}
\def\goto{{\rightarrow}}
 \def\ar{\!\!&}\def\nnm{\nonumber}\def\ccr{\nnm\\}

\begin{document}

\title{\large \bf Well-posedness of martingale problem for SBM\\
 with interacting branching}

\author{Lina Ji\footnote{Department of Mathematics, Southern University of Science and Technology, Shenzhen, China. Email: jiln@sustech.edu.cn.}, Jie Xiong\footnote{Department of Mathematics \& National Center for Applied Mathematics (Shenzhen), Southern University of Science and Technology, Shenzhen, China. Email: xiongj@sustech.edu.cn.} and Xu Yang\footnote{School of Mathematics and Information Science, North Minzu University, Yinchuan, China. Email: xuyang@mail.bnu.edu.cn.}}

\date{\today}
\maketitle
\begin{abstract}

In this paper a martingale problem for super-Brownian motion with interactive branching is derived. The uniqueness of the solution to the martingale problem is obtained by using the pathwise uniqueness of the solution to a corresponding system of SPDEs with proper boundary conditions. The existence of the solution to the martingale problem and the H\"{o}lder continuity of the density process are also studied.

\end{abstract}

\bf Keywords. \rm super-Brownian motion; interacting branching; function-valued process; stochastic partial differential equation.

\bf AMS Mathematics subject classification (2010): \rm 60H15; 60J68.

\section{Introduction and main results}
\setcounter{equation}{0}
\renewcommand{\theequation}{\thesection.\arabic{equation}}

Let $M(\mbb{R})$ be the collection of all finite Borel measures on $\mbb{R}.$ Let $C_b^n(\mbb{R})$ be the collection of all bounded continuous functions on $\mbb{R}$ with bounded derivatives up to $n$th order.
We consider a  $M(\mbb{R})$-valued process $(X_t)_{t \ge 0}$ satisfying the following martingale problem (MP): for $\phi \in C_b^2(\mbb{R}),$ the process
\beqlb\label{0.1}
M_t(\phi)\equiv\<X_t,\phi\> - \<X_0, \phi\>-\int^t_0\<X_s,\frac12\phi''\>ds
\eeqlb
is a continuous martingale with
\beqlb\label{0.2}
\<M(\phi)\>_t=\int_0^t\<X_s,\ga(\mu_s,\cdot)\phi^2\>ds,
\eeqlb
where $\gamma$ is the interacting branching rate depending on the density process. The notation $\< \nu, \phi\>$ denotes the integral of the function $\phi$ with respect to the measure $\nu.$ In this paper we assume $X_0(d x) = \mu_0(x)d x$ with $\mu_0 \in L^2(\mbb{R})$, where $L^2(\mbb{R})$ is the space of functions $f$ such that $\int_\mbb{R}f(x)^2 d x < \infty.$ 
When $\gamma$ is a constant, the process $(X_t)_{t \ge 0}$ is a super-Brownian motion. In this case the well-posedness of the MP (\ref{0.1}, \ref{0.2}) was established by the nonlinear partial differential equation satisfied by its log-Laplace transform. Moreover, a new approach for the well-posedness of the MP was suggested by Xiong~\cite{X13b},  in which a relationship between the super-Brownian motion and a  stochastic partial differential equation (SPDE) satisfied by its corresponding distribution-function-valued process was established. The weak uniqueness of the solution to the MP~(\ref{0.1}, \ref{0.2}) was also obtained by the strong uniqueness of the solution to the corresponding SPDE in  \cite{X13b}. See He et.~al~\cite{HLY14} for the case of super-L\'{e}vy process. 


While the superprocesses with interaction are more natural since the branching and spatial motion for many species  depend on the processes themselves. When the spatial motion is interactive,  the well-posedness of the martingale problem was studied by Donnelley and Kurtz~\cite{DK99}, see also 
 Perkins~\cite[Theorem~V.5.1]{P95} and Li et.~al \cite{LWX05}. Uniqueness for the historical superprocesses with certain interaction was investigated by Perkins~\cite{P95}. Furthermore, the superprocesses with interactive immigration were studied by \cite{ DL03, FL04,S90}, see also Li~\cite[Section~10]{L11}. The well-posedness of the martingale problem for the interactive immigration process was solved by Mytnik and Xiong \cite{MX15}. See also \cite{XY16} for the well-posedness of the martingale problem for a superprocess with location-dependent branching, interactive immigration mechanism and spatial motion.

However, the hard case for the superprocess with interactive branching was rarely investigated.
 We are interested in the case of $\gamma(\mu_s, x) = \gamma(\mu_s(x)),$ i.e., the interactive branching mechanism depending on the density process $(\mu_t(x))_{t \ge 0, x \in \mbb{R}}.$ However, in this case the well-posedness of the MP (\ref{0.1}, \ref{0.2}) is still an {\it open problem}. The weak uniqueness of the solution to the MP (\ref{0.1}, \ref{0.2}) is very difficult to prove. As a first step, throughout this paper we assume $\gamma$ satisfies the following condition:

\begin{con}\label{0.4b}
Fixed integer $n \ge 0.$ For $-\infty=a_0<a_1<\cdots<a_n<a_{n+1}=\infty$, let 
 \beqnn
\ga(\mu_s,x)=
\begin{cases}
g_i^2(\mu_s(a_{i + 1})),& a_i\le x < a_{i+1},\ i = 0, \cdots, n - 1,\\
 g_n^2\(\int_{a_n}^\infty\mu_s(y)d y\), & a_n \le x < \infty,
\end{cases}
\eeqnn
where $g_i,\ i = 0, \cdots, n$ are positive continuous bounded functions from $\mbb{R}_+$ to $\mbb{R}_+$.
\end{con}
The existence and H\"{o}lder continuity of the density process $(\mu_t(x))_{t \ge 0, x \in \mbb{R}}$ are investigated in this paper, which satisfies the following SPDE:
\beqlb\label{1.1}
\frac{\partial}{\partial t} \mu_t(x) = \frac{1}{2}\Delta \mu_t(x) + \sqrt{\mu_t(x)\gamma(\mu_t, x)}\dot{W}(t, x),
\eeqlb
where $\Delta$ denotes the one-dimensional Laplacian operator and $\{W(t, x): t \ge 0, x \in \mbb{R}\}$ is a time-space Gaussian white noise based on the Lebesgue measure and the dot denotes the derivative in distribution sense. The existence of the solution to the MP (\ref{0.1}, \ref{0.2}) is given by showing the existence of the solution to \eqref{1.1}. Moreover, we show the weak uniqueness of the solution to the MP~(\ref{0.1}, \ref{0.2}) in Theorem~\ref{main}. The main idea is to relate the MP with a system of SPDEs, which is satisfied by a sequence of corresponding function-valued processes on intervals. The weak uniqueness of solution to the MP follows from the pathwise uniqueness of the solution to the system of SPDEs, see Section 3.  Throughout this paper we always assume that all random variables defined on the same filtered probability space $(\Omega, \mcr{F}, \mcr{F}_t, \mbb{P}).$ Let $\mbb{E}$ denote the corresponding expectation.

The main results of this paper without special explanation are stated as below.

\begin{thm}\label{t1} (Existence)
There exists a $M(\mbb{R})$-valued continuous process $(X_t)_{t \ge 0}$ satisfying the MP (\ref{0.1}, \ref{0.2}). Moreover, $X_t(d x)$ is absolutely continuous respect to $d x$ with density $\mu_t(x)$ satisfying \eqref{1.1}.
\end{thm}

\begin{thm}\label{t2}
(Joint H{\"o}lder continuity)
Suppose that $(\mu_t(x))_{t \ge 0, x \in \mbb{R}}$ satisfies \eqref{1.1} and $\mu_0$ is H{\"o}lder continuous with exponent $0 < \lambda < \frac{1}{2}$. Then $[0,  T] \times \mbb{R} \ni (t, x) \goto \mu_t(x)$ is  H{\"o}lder continuous with exponent $\lambda_1/2$ in time variable and with exponent $\lambda_2$ in space variable, where $\lambda_1, \lambda_2 \in (0, 1/2).$ Namely, there exists a random variable $K \ge 0$ only depending on $\lambda_1$ and $\lambda_2$ such that
\beqnn
|\mu_t(x) - \mu_r(y)| \le K (|t - r|^{\lambda_1/2} + |x - y|^{\lambda_2}), \qquad t, r \in [0, T],\ x, y \in \mbb{R}.
\eeqnn
\end{thm}

\begin{thm}\label{main} (Uniqueness)
Assume $g_n$ is $\beta$-H\"{o}lder continuous with $ \frac{1}{2}\le \beta \le 1$, i.e., 
\beqlb\label{beta}
|g_n(x) - g_n(y)| \le K |x - y|^\beta, \qquad  x, y \ge 0
\eeqlb
for some constant $K.$ Then the weak uniqueness of the solution to the MP (\ref{0.1}, \ref{0.2}) holds.
\end{thm}

\begin{rem}
Taking $n = 0$ in Condition~\ref{0.4b}, the branching rate depends on the total mass process, i.e., the quadratic variation process of the martingale defined by \eqref{0.1} is
\beqlb\label{rem}
\langle M(\phi)\rangle_t = \int_0^t \gamma(\langle X_s, 1\rangle )\langle X_s, \phi^2\rangle d s.
\eeqlb
The well-posedness of the MP (\ref{0.1}, \ref{rem}) is a corollary of Theorem~\ref{t1} and \ref{main}.
\end{rem}

We introduce some notation. Let $\mathcal{X}_0$ be the Hilbert space consisting of all functions $f$ such that
\beqnn
\|f\|_0^2 :=  \int_{\mbb{R}}f(x)^2e^{-|x|}d x < \infty.
\eeqnn
We denote the corresponding inner product by $\<\cdot, \cdot\>_0.$
Let $\mathscr{B}(\mbb{R})$ (resp. $\mathscr{B}(a_1, a_2)$) denote the Borel $\sigma$-algebra on $\mbb{R}$ (resp. $(a_1, a_2)$). Let $B(\mbb{R})^+$ (resp. $C_b(\mbb{R})^+$) be the collection of all bounded positive (resp. bounded positive continuous) functions on $\mbb{R}.$
Let $B[a_1, a_2]$ be the Banach space of bounded measurable functions on $[a_1, a_2]$. For $f, g \in B [a_1, a_2]$ let $\< f, g\> = \int_{a_1}^{a_2} f(x)g(x)\d x$. Define $C_b[a_1, a_2]$ be the set of bounded continuous functions on $[a_1, a_2]$. For any integer $n \ge 0,$ let $C_b^n[a_1, a_2]$ be the subset of $C_b[a_1, a_2]$ of functions with bounded continuous derivatives up to the $n$th order. Let $C_c^n(a_1, a_2)$ denote the subset of $C_b^n[a_1, a_2]$ of functions with compact supports in $(a_1, a_2).$


The rest of the paper is organized as follows. In Section~2, we give the proofs of Theorem~\ref{t1} and \ref{t2}. The weak uniqueness of the solution to the MP~(\ref{0.1}, \ref{0.2}), i.e., the proof of Theorem~\ref{main}, is given in Section~3. Throughout the paper we use $\nabla$ to be the first order spatial differential operator, use $K$ to denote a non-negative constant whose value may change from line to line. In the integrals, we make convention that, for $a \le b \in \mbb{R},$
\beqnn
\int_a^b = \int_{(a, b]}\quad \text{and}\quad \int_a^\infty = \int_{(a, \infty)}.
\eeqnn

\section{Proofs of Theorem~\ref{t1} and \ref{t2}}

In this section, we give the proofs of Theorem~\ref{t1} and \ref{t2}. The existence of the solution to MP (\ref{0.1}, \ref{0.2}) is obtained by the existence of the corresponding density process. Moreover, the H{\"o}lder continuity of density process $(\mu_t(x))_{t \ge 0, x \in \mbb{R}}$ is given by a standard argument.



\begin{lemma}\label{t2.3}
The martingale defined in (\ref{0.1},\ref{0.2}) induces a $(\mcr{F}_t)$-martingale measure $\{M_t(B): t \ge 0, B \in \mcr{B}(\mbb{R})\}$ satisfying
\beqnn
M_t(\phi) = \int_0^t \int_{\mbb{R}} \phi(x) M(d s, d x), \qquad t \ge 0,\ \phi \in C_b^2(\mbb{R}),
\eeqnn
where $M(d s, d x)$ is an orthogonal martingale measure on $\mbb{R}_+\times\mbb{R}$ with covariance measure $d s\int_{\mbb{R}}\big[\gamma(X_s, z) \delta_z(d x)\delta_z(d y)\big]X_s(d z)$.
\end{lemma}

\proof Notice that $\gamma$ is bounded by  Condition \ref{0.4b}. By the MP~(\ref{0.1}, \ref{0.2}), one can see that $\mbb{E}[\<X_t, 1\>] = \<X_0, 1\> < \infty.$ For each $n \ge 1$ we define the measure $\Gamma_n \in M(\mbb{R})$ by 
\beqnn
\Gamma_n(\phi) = \mbb{E}\Big[\int_0^n d s \int_{\mbb{R}}\gamma(\mu_s, x) \phi(x)X_s(d x)\Big]
\eeqnn
with $\phi \in B(\mbb{R})^+.$ Then $\Gamma_n(\phi) \le K \int_0^n \mbb{E} [\<X_s, 1\>] d s$ is bounded for each $n \ge 1.$ The rest proof follows by changing $c(z) = \gamma(\mu_s, z)$ and $H(z, \d \nu) = 0$ in the proof of Li~\cite[Theorem~7.25]{L11}. We omit it here.
\qed

Let $T_t(x, d y)$ be the semigroup generated by $\Delta/2,$ which is absolutely continuous with respect to Lebesgue measure $\d y$ with density $p_t(x, y)$ satisfying
\[p_t(x, y) = p_t(x - y) =\frac{1}{\sqrt{2\pi t}}\e^{-|x - y|^2/(2t)}.\]
 
\begin{lemma}\label{<mu,1>}
 Suppose that $(\mu_t(x))_{t \ge 0, x \in \mbb{R}}$ is a solution to \eqref{1.1} with $\mu_0 \in L^2(\mbb{R})$. Then for every $t \ge 0,$ we have $\mbb{E}[\langle\mu_t, 1\rangle] < \infty.$
 \end{lemma}
 
 \proof
 By \eqref{1.1} we have
\beqnn
\langle \mu_t, 1\rangle = \langle \mu_0, 1\rangle + \int_0^t\int_{\mbb{R}} \sqrt{\mu_s(x)\gamma(\mu_s, x)} W(d s, d x).
\eeqnn 
Notice that  $\sqrt{x} \le x + 1$ for any $x \ge 0.$ By H{\"o}lder's inequality one can check that
\beqnn
\mbb{E}[\langle \mu_t, 1\rangle] \ar\le\ar \langle \mu_0, 1\rangle + \mbb{E}\left|\int_0^t\int_{\mbb{R}} \sqrt{\mu_s(x)\gamma(\mu_s, x)} W(d s, d x)\right|\cr
\ar\le\ar \langle \mu_0, 1\rangle + K\left[\mbb{E}\left(\int_0^t\int_{\mbb{R}} \mu_s(x) d s d x\right)\right]^{1/2}\cr
\ar\le\ar K + K\int_0^t \mbb{E}\left[\langle \mu_s, 1\rangle \right]d s.
\eeqnn
The result follows by Gronwall's inequality.
 \qed  
 
 Similar with Li~\cite[Theorem~7.26, Theorem~7.28]{L11}, we have the following result.
 
\begin{prop}\label{t2.4}
Suppose that $(X_t)_{t \ge 0}$ is a solution to the MP (\ref{0.1}, \ref{0.2}). Then for each $t \ge 0$ the random measure $X_t(d x)$ is absolutely continuous respect to $d x$ with density $\mu_t(x)$ satisfying \eqref{1.1}.
Conversely, assume that $(\mu_t(x))_{t \ge 0, x \in \mbb{R}}$ is a solution to \eqref{1.1} with $\mu_0 \in L^2(\mbb{R})$. Then $(X_t)_{t \ge 0}$ satisfies the MP (\ref{0.1}, \ref{0.2}).
\end{prop}

\proof
Suppose that $(X_t)_{t \ge 0}$ satisfies the MP (\ref{0.1}, \ref{0.2}). By Lemma~\ref{t2.3} and the proof of \cite[Theorem~7.26]{L11}, we have
\beqnn
\<X_t, \phi\> = \<X_0, T_t\phi\> + \int_0^t\int_{\mbb{R}}T_{t - s}\phi(x)M(d s, d x),
\eeqnn
where $M(d s, d x)$ is an orthogonal martingale measure defined in Lemma~\ref{t2.3}. Recall that $\gamma$ is bounded by Condition~\ref{0.4b}. Let $\phi \in C_b(\mbb{R})^+$ and set $p_t (x, z) = 0$ for all $t \le 0.$ For any $n \ge 1$ we have
\beqnn
&\ &\mbb{E}\Big[\int_{\mbb{R}}\phi(z)d z\int_0^n \int_{\mbb{R}}p_{t - s}(x - z)^2 {\gamma(\mu_s, x)}X_s(d x)d s\Big]\\
&\ &\qquad \le  K \mbb{E}\Big[\int_{\mbb{R}}\phi(z)d z\int_0^t \int_{\mbb{R}}p_{t - s}(x - z)^2 X_s(d x)d s\Big]\\
&\ &\qquad \le K \sup_{z \in \mbb{R}}\phi(z)\int_0^t \frac{\mbb{E}\left[\<X_s, 1\>\right]}{\sqrt{2\pi (t - s)}} \d s = K \sup_{z \in \mbb{R}}\phi(z) \sqrt{t} < \infty.
\eeqnn
Then by stochastic Fubini's theorem (e.g., see Li~\cite[Theorem~7.24]{L11}), we get $\<X_t, \phi\> = \int_\mbb{R} \mu_t(x)\phi(x)d x$
with
\beqnn
\mu_t(x) = \int_\mbb{R}p_t(x - z)\mu_0(z)\d z + \int_0^t\int_{\mbb{R}}p_{t - s}(x - z)M(d s, d z).
\eeqnn
By El Karoui and M\'{e}l\'{e}ard \cite[Theorem III-6]{EM90}, on some extension of the probability space one can define a white noise $W(d s, d z)$ on $\mbb{R}_+\times \mbb{R}$ based on $d s d z$ such that
the following holds:
\beqlb\label{mu}
\mu_t(x) = \<\mu_0, p_t(x - \cdot)\>  + \int_0^t\int_{\mbb{R}}\sqrt{\mu_s(z)\gamma(\mu_s, z)}p_{t - s}(x - z)W(d s, d z),
\eeqlb
 which implies \eqref{1.1}.

Conversely, suppose that $\mu_t(x)$ satisfies \eqref{1.1} with $\mu_0 \in L^2(\mbb{R})$ and denote $X_t(d x) = \mu_t(x)d x$. It follows from Lemma~\ref{<mu,1>} that $X_t \in M(\mbb{R})$ almost surely for every $t \ge 0$. For any $\phi \in C_b^2(\mbb{R})$ one can check that
\beqnn
\<X_t, \phi\> \ar=\ar \int_{\mbb{R}}\mu_t(x)\phi(x) d x\\
 \ar=\ar \int_{\mbb{R}}\mu_0(x)\phi(x) d x + \frac{1}{2}\int_0^t \int_{\mbb{R}}\mu_s(x)\phi''(x) d x d s\\
\ar\ar +\int_0^t\int_{\mbb{R}}\phi(x)\sqrt{\mu_s(x)\gamma(\mu_s, x)}W(d s, d x)\\
\ar=\ar \<X_0, \phi\> + \frac12\int_0^t\<X_s, \phi''\>d s + M_t(\phi)
\eeqnn
with $\<M(\phi)\>_t$ satisfying \eqref{0.2}, which completes the proof.
\qed

Now we show the existence of the solution to \eqref{1.1}. For any $T > 0,$ let $m \ge 1$ and $t_k = kT/m$ with $k = 0, 1, \cdots, m.$ Define a sequence of approximation by
 \beqlb\label{mun}
\mu_t^{m}(x) = \mu_0(x) + \frac12\int_0^t\Delta \mu_s^m(x) d s + \sum_{k = 1}^m\sqrt{\gamma(\mu_{t_{k - 1}}^m, x)}\int_{ t_{k - 1}\wedge t}^{ t_k\wedge t}  G_m(\mu_s^m(x))\dot{W}(s, x)d s,
\eeqlb
where
\beqnn
G_m(x) = \int_{\mbb{R}} p_{m^{-1}}(x - y)(\sqrt{|y|}\wedge m)d y 
\eeqnn
is a Lipschitz function on $[0, \infty)$ for fixed $m \ge 1$ and $\lim_{m \rightarrow \infty}G_m(x) = \sqrt{x}$ for all $x \ge 0.$ For all $m \ge 1$ and $x \ge 0,$ one can check that there is a constant $K > 0$ such that 
\beqlb\label{G_m}
G_m(x) \ar\le\ar \int_{\mbb{R}}p_{m^{-1}}(x - y)(1 + |y|)d y = 1 + \mbb{E}[|B_{m^{-1}}^x|] \le 1 + {  \left[\mbb{E}\left(|B_{m^{-1}}^x|^2\right)\right]^{1/2}}\cr
\ar\le\ar 1 + (1 + x^2)^{1/2} \le K(x + 1)
\eeqlb
 where $(B_t^x)_{t \ge 0}$ is a Brownian motion with initial value $x.$ For any $k = 1, \cdots, m,$ conditioned on $\mcr{F}_{t_{k - 1}},$ it is well known that there is a strong unique non-negative solution $\{\mu_t^m(x): t \in [t_{k - 1}, t_k), x \in \mbb{R}\}$ to
 \beqnn
\mu_t^{m}(x) = \mu_{t_{k - 1}}^m(x) + \frac12\int_{t_{k - 1}}^t\Delta \mu_s^m(x) d s + \sqrt{\gamma(\mu_{t_{k - 1}}^m, x)}\int_{t_{k - 1}}^t G_m(\mu_s^m(x))
\dot{W}(s, x)d s,
\eeqnn
 see \cite[Theorem~5.1]{D75}. Therefore, the process $(\mu_t^{m}(x))_{t \ge 0, x \in \mbb{R}}$ is the unique strong solution to \eqref{mun}. 

Let $J(x) = \int_{\mbb{R}}e^{-|y|}\rho(x - y) d y,$
where $\rho$ is the mollifier given by $\rho(x) = C\exp\{-1/(1 - x^2)\}1_{\{|x| < 1\}},$ and $C$ is a constant  such that $\int_{\mbb{R}}\rho(x) d x = 1.$ Then for any $n \ge 1,$ there are constants $c_n$ and $C_n$ such that
\beqlb\label{J}
c_n e^{-|x|} \le J^{(n)}(x) \le C_n e^{-|x|}, \qquad \forall\ x \in \mbb{R},
\eeqlb
see (2.1) of Mitoma~\cite{M85}. We may and will replace $e^{-|x|}$ by $J(x)$ in the definition of the space $\mathcal{X}_0,$ and define $\<f, g\>_0 = \int_{\mbb{R}}f(x)g(x)J(x)d x$ for any $f, g \in C_b^2(\mbb{R})$.

\begin{lemma}\label{l2.4}
Assume that $\mu_0$ satisfies $\int_{\mbb{R}}\mu_0(x)^{2p}J(x)dx < \infty$ for some $p \ge 1$. Then for every $T > 0$, we have
\beqnn
\sup_{0 \le t \le T, m \ge 1}\mbb{E}\left[\int_{\mbb{R}}\mu_t^m(x)^{2p}J(x)d x\right] < \infty.
\eeqnn
\end{lemma}

\proof
Using the convolution form, the solution $\mu_t^m(x)$ to \eqref{mun} can be represented as
\beqlb\label{mun1}
\mu_t^m(x) \ar=\ar \<\mu_0, p_t(x - \cdot)\>\ccr
\ar\ar + \sum_{k = 1}^m\int_{t_{k - 1}\wedge t}^{t_k\wedge t}\int_{\mbb{R}} \sqrt{\gamma(\mu_{t_{k - 1}}^m, z)}G_m(\mu_s^m(z))p_{t - s}(x - z) W(d s, d z).
\eeqlb
By H{\"o}lder's inequality we have
\beqnn
\ar\ar\left|\int_0^t d s\int_{\mbb{R}} \mu_s^m(z)^2p_{t - s}(x - z)^2 d z\right|^p\\
\ar\ar\qquad\le \left|\int_0^t \frac{1}{\sqrt{2\pi(t - s)}}d s\int_{\mbb{R}} \mu_s^m(z)^2p_{t - s}(x - z) d z\right|^p\\
\ar\ar \qquad\le K\int_0^t \frac{1}{\sqrt{2\pi(t - s)}}d s\left|\int_{\mbb{R}} \mu_s^m(z)^{2}p_{t - s}(x - z) d z\right|^p \cdot \left|\int_0^t \frac{1}{\sqrt{2\pi(t - s)}}d s\right|^{p/q}\\
\ar\ar \qquad\le K\int_0^t \frac{1}{\sqrt{2\pi(t - s)}}d s\int_{\mbb{R}} \mu_s^m(z)^{2p}p_{t - s}(x - z) d z 
\eeqnn
with $1/p + 1/q = 1$ and $p, q \ge 1.$  
By the above inequality, \eqref{G_m} and Burkholder-Davis-Gundy's inequality one can see that
\beqnn
\ar\ar\mbb{E}\left[\int_{\mbb{R}}J(x)d x\left|\sum_{k = 1}^m\int_{t_{k - 1}\wedge t}^{t_k\wedge t}\int_{\mbb{R}} \sqrt{\gamma(\mu_{t_{k - 1}}^m, z)}G_m(\mu_s^m(z))p_{t - s}(x - z) W(d s, d z)\right|^{2p}\right]\\
\ar\ar\qquad \le K \mbb{E}\left[\int_{\mbb{R}}J(x)d x\left|\int_0^t \int_{\mbb{R}} G_m(\mu_s^m(z))^2p_{t - s}(x - z)^2 d s d z\right|^p\right]\\
\ar\ar\qquad \le K \mbb{E}\left[\int_{\mbb{R}}J(x)d x\left|\int_0^t d s \int_{\mbb{R}} (\mu_s^m(z)^2 + 1)p_{t - s}(x - z)^2 d z\right|^p\right]\\
\ar\ar\qquad { \le K +  K \mbb{E}\left[\int_0^t \frac{1}{\sqrt{2\pi (t - s)}} d s \int_{\mbb{R}}J(x)d x\int_{\mbb{R}} \mu_s^m(z)^{2p}p_{t - s}(x - z) d z\right]}\\
\ar\ar\qquad { \le K +  K \mbb{E}\left[\int_0^t \frac{1}{\sqrt{2\pi (t - s)}} d s \int_{\mbb{R}}\mu_s^m(z)^{2p} d z\int_{\mbb{R}}J(u + z) p_{t - s}(u) d u\right]}\\
\ar\ar\qquad { \le K +  K \int_0^t \frac{1}{\sqrt{t - s}}\mbb{E}\left[\int_{\mbb{R}}\mu_s^m(z)^{2p} J(z) d z\right] d s.}
\eeqnn
By \eqref{J} we have $J(u + z) \le J(z)e^{|u|}$ which implies the last inequality, since
\beqnn
\int_{\mbb{R}}e^{u}p_{t - s}(u)d u \le \int_{\mbb{R}}e^{T|u|}p_1(u)d u = \mbb{E}\left[e^{T|B_1|}\right] \le K
\eeqnn
for $0 \le s \le t \le T,$ where $B_t$ is a standard Brownian motion. Moreover, we have
\beqnn
\int_{\mbb{R}}\<\mu_0, p_t(x - \cdot)\>^{2p}J(x)d x \ar=\ar \int_{\mbb{R}}J(x)d x\left[\int_{\mbb{R}}p_t(x - z)\mu_0(z)d z\right]^{2p}\\
\ar\le\ar \int_{\mbb{R}}J(x)d x\int_{\mbb{R}}p_t(x - z)\mu_0(z)^{2p}d z\\
\ar\le\ar K\int_{\mbb{R}}\mu_0(z)^{2p}J(z)d z < \infty.
\eeqnn
Combining the above inequality, we have
\beqnn
\mbb{E}\left[\int_{\mbb{R}}\mu_t^m(x)^{2p}J(x)d x\right] \le K + K\int_0^t \frac{1}{\sqrt{t - s}} \mbb{E}\left[\int_{\mbb{R}}\mu_s^m(x)^{2p}J(x)d x\right]d s.
\eeqnn
Iterating the above once, one can check that
{\small\beqnn
\mbb{E}\left[\int_{\mbb{R}}\mu_t^m(x)^{2p}J(x)d x\right] \ar\le\ar K + K\int_0^t \mbb{E}\left[\int_{\mbb{R}}\mu_r^m(x)^{2p}J(x)d x\right]d r\int_r^t \frac{1}{\sqrt{(t - s)(s - r)}}d s\\
\ar\le\ar  K + K\int_0^t \mbb{E}\left[\int_{\mbb{R}}\mu_r^m(x)^{2p}J(x)d x\right]d r.
\eeqnn}
The result follows by Gronwall's inequality. 
\qed



We proceed to proving the tightness of $(\mu^m_\cdot)$ in $C([0, T] \times \mbb{R}).$ Denote
\beqlb\label{num}
\nu_t^m(x) = \sum_{k = 1}^m\int_{t_{k - 1}\wedge t}^{t_k\wedge t}\int_{\mbb{R}} \sqrt{\gamma(\mu_{t_{k - 1}}^m, z)}G_m(\mu_s^m(z))p_{t - s}(x - z) W(d s, d z).
\eeqlb

\begin{lemma}\label{l22.6}
Assume that $\mu_0$ satisfies $\int_{\mbb{R}}\mu_0(x)^{2p}J(x)dx < \infty$ for some $p \ge 1$. For fixed $0 < \alpha < 1$ and $T > 0,$ there is a constant $K$ such that
\beqnn
\mbb{E}[|\nu_t^m(x) - \nu_r^m(x)|^{2p}]\le K(t - r)^{\alpha p/2}
\eeqnn
for all $0 < r < t \le T.$
\end{lemma}
\proof
Recall that $\gamma$ is bounded by Condition \ref{0.4b}. By \eqref{G_m}, \eqref{num} and Burkholder-Davis-Gundy's inequality one can check that
\beqnn
\mbb{E}\left[|\nu_t^m(x) - \nu^m_r(x)|^{2p}\right]\ar\le\ar  K \mbb{E}\left[\left|\int_0^r\int_{\mbb{R}}\left[p_{r - s}(x - z) - p_{t - s}(x - z)\right]^2\left[\mu_s^m(z)+1\right] d s d z\right|^p\right]\\
\ar\ar + K\mbb{E}\left[\Big|\int_r^t\int_{\mbb{R}}[\mu_s^m(z) + 1] p_{t - s}(x - z)^2d s d z\Big|^p\right]\\
\ar=\ar I_1^m(t,r) + I_2^m(t,r).
\eeqnn
By \cite[Lemma~III4.5]{P02}, for any $0 < s < r < t \le T$ and $\delta \in (0, 1/4),$ there is a constant $K = K(T) > 0$ such that
\beqlb\label{p_{r - s}}
\ar\ar\left[p_{r - s}(x - z) - p_{t - s}(x - z)\right]^2\cr
\ar\ar\qquad \le K (t - r)^{\delta}(r - s)^{-3\delta/2}\left[p_{r - s}(x - z)^{2 - \delta} + p_{t - s}(x - z)^{2 - \delta}\right].
\eeqlb 
By \cite[Lemma~1.4.4]{X13b}, we have
\beqnn
\int_0^r\int_{\mbb{R}}[p_{r - s}(x - z) - p_{t - s}(x - z)]^2 d s d z \le K|t - r|^{1/2}
\eeqnn
and $\int_r^t\int_{\mbb{R}} p_{t - s}(x - z)^2d s d z \le K|t - r|^{1/2}.$ From \eqref{p_{r - s}} and Lemma~\ref{l2.4} one can obtain that
\beqnn
&\ &\mbb{E}\left[\int_0^rd s\int_{\mbb{R}}[p_{r - s}(x - z) - p_{t - s}(x - z)]^2|\mu_s^m(z)|^pd z\right] \\
&\ &\quad\le K(t - r)^{\delta}\int_0^r(r - s)^{-3\delta/2}d s\mbb{E}\left[\int_{\mbb{R}}[p_{r - s}(x - z)^{2 - \delta} + p_{t - s}(x - z)^{2 - \delta}]|\mu_s^m(z)|^p d z\right]\\
&\ &\quad\le K(t - r)^{\delta}\int_0^r(r - s)^{-3\delta/2}d s\left(\int_{\mbb{R}}[p_{r - s}(x - z)^{4 - 2\delta} + p_{t - s}(x - z)^{4 - 2\delta}]e^{|z|}d z\right)^{1/2}\\
&\ &\qquad \cdot\mbb{E}\left[\left(\int_{\mbb{R}}|\mu_s^m(z)|^{2p}J(z) d z\right)^{1/2}\right]\\
&\ &\quad \le K (t - r)^{\delta}\left[\int_0^r (r - s)^{-(3 + 4\delta)/4}d s + \int_0^r (r - s)^{-3\delta/2}(t - s)^{-(3 - 2\delta)/4}d s\right]\\
&\ &\quad \le  K (t - r)^{\delta}.
\eeqnn
By H\"{o}lder's inequality, it implies that
\beqlb\label{I_1}
I_1^m(t,r) \ar\le\ar K\mbb{E}\left[\int_0^rd s\int_{\mbb{R}}[p_{r - s}(x - z) - p_{t - s}(x - z)]^2|\mu_s^m(z)|^pd z\right]\cr
\ar\ar\quad \cdot \left[\int_0^r\int_{\mbb{R}}[p_{r - s}(x - z) - p_{t - s}(x - z)]^2 d s d z\right]^{p - 1}
 + K(t - r)^{p/2}\ccr
\ar\le\ar K (t - r)^{\delta + (p-1)/2} + K(t - r)^{p/2}.
\eeqlb 
Similarly, by H\"{o}lder's inequality and Lemma~\ref{l2.4} one can see that
\beqnn
\ar\ar\mbb{E}\left[\int_r^t \int_{\mbb{R}}|\mu_s^m(z)|^p p_{t - s}(x - z)^2 d s d z\right]\ccr
\ar\ar\qquad \le \mbb{E}\left[\int_r^t d s\left[\int_{\mbb{R}}|\mu_s^m(z)|^{2p} J(z)d z\right]^{1/2}\left[\int_{\mbb{R}}p_{t - s}(x - z)^4 e^{|z|} d z\right]^{1/2}\right]\ccr
\ar\ar\qquad \le K\int_r^t (t - s)^{-3/4}d s = K(t - r)^{1/4},
\eeqnn
which implies that
\beqlb\label{I_2}
I_2^m(t,r) \ar\le\ar K \mbb{E}\left[\int_r^t \int_{\mbb{R}}|\mu_s^m(z)|^p p_{t - s}(x - z)^2 d s d z\right]\left[\int_r^t \int_{\mbb{R}} p_{t - s}(x - z)^2 d s d z\right]^{p - 1}\ccr
\ar\ar + K(t - r)^{p/2}\ccr
\ar\le\ar K(t - r)^{1/4 + (p - 1)/2} + K(t - r)^{p/2}.
\eeqlb 
The result follows.
\qed

Similar with above, we have the following result.

\begin{lemma}\label{l22.7}
Assume that $\mu_0$ satisfies $\int_{\mbb{R}}\mu_0(x)^{2p}J(x)dx < \infty$ for some $p \ge 1$. For fixed $0 < \beta < 1,$ there is a constant $K$ such that
\beqnn
\mbb{E}\left[|\nu_t^m(x) - \nu_t^m(y)|^{2p}\right]\le K|x - y|^{\beta p}
\eeqnn
for any $x, y \in \mbb{R}$.
\end{lemma}


{\bf Proof of Theorem \ref{t1}.} By \eqref{mun1} one can check that
\beqnn
\mu_t^m(x) \ar=\ar \<\mu_0, p_t(x - \cdot)\> + \nu_t^m(x).
\eeqnn
By Kolmogorov's criteria (see \cite[Corollary 16.9]{K02}), Lemma~\ref{l22.6} and Lemma~\ref{l22.7}, for each fixed $T, K > 0,$ the sequence of laws of $\{\mu_t^m(x): (t, x) \in [0, T]\times [-K, K]\}$ on $C([0, T] \times [-K, K])$ is tight, and hence, has a convergent subsequence. By the standard diagonlization argument, there exists a subsequence $(\mu_t^{m_k}, W_t^{m_k})$ which converges to  $(\mu_t, W_t)$ in law on $C([0, T]\times[-K,K])$ for each $K$ and $T.$ Therefore, $(\mu_t^{m_k}, W_t^{n_k})_{t \ge 0}$ converges in law as $k \rightarrow \infty.$ Applying Skorokhod's representation (e.g. \cite[Theorem~1.8]{EK86}), on another probability space, there are continuous processes $(\hat{\mu}_t^{m_k}, \hat{W}_t^{m_k})_{t \ge 0}$ and $(\hat{\mu}_t, \hat{W}_t)_{t \ge 0}$
with the same distribution as  $(\mu_t^{m_k}, W_t^{m_k})$ and  $(\mu_t, W_t),$ respectively. Moreover,
\beqnn
(\hat{\mu}_t^{m_k}, \hat{W}_t^{m_k})_{t \ge 0} \rightarrow (\hat{\mu}_t, \hat{W}_t)_{t \ge 0}
\eeqnn
 almost surely as $k \rightarrow \infty$. Recall that there is a unique strong non-negative solution to \eqref{mun}. For any $f \in C_b^2(\mbb{R})$ with $|f| + |f''| \le K J$, we have
 \beqnn
 \<\hat{\mu}_t^{m_k}, f\> \ar=\ar \<\mu_0, f\> + \frac12\int_0^t \<\hat{\mu}_s^{m_k}, f''\>d s \\
 &\ &  + \sum_{j = 1}^{m_k}\int_{t_{j - 1}\wedge t}^{t_j\wedge t}\int_{\mbb{R}} \sqrt{\gamma(\hat{\mu}_{t_{j - 1}}^{m_k}, z)}G_{m_k}(\hat{\mu}_s^{m_k}(z))f(z) \hat{W}^{m_k}(d s, d z).
 \eeqnn
It follows from Lemma~\ref{l2.4} that for every $t > 0$, we have
\beqlb\label{J1}
\lim_{k \rightarrow \infty}\mbb{E}\left[\int_{\mbb{R}}|\hat{\mu}_t^{m_k}(x) - \hat{\mu}_t(x)|^{2}J(x)d x\right] = 0
\eeqlb
and
\beqlb\label{J2}
\sup_{0 \le t \le T, k \ge 1}\mbb{E}\left[\int_{\mbb{R}}[\hat{\mu}_t(x) + \hat{\mu}_t^{m_k}(x)]^{2}J(x)d x\right] < \infty,
\eeqlb
since $\mu_0(x) \in L^2(\mbb{R}).$ Then by the above and dominated convergence theorem,  we have
\beqnn
\mbb{E}\left[\int_0^t\langle|\hat{\mu}_s^{m_k} - \hat{\mu}_s|, f''\rangle d s\right] \le \mbb{E}\left[\int_0^td s\int_{\mbb{R}}|\hat{\mu}_s^{m_k}(x) - \hat{\mu}_s(x)|^2J(x)d x\right]^{1/2} \rightarrow 0
\eeqnn
as $k \rightarrow \infty.$ Recall that $\gamma$ satisfies Condition \ref{0.4b}. Thus, 
\beqnn
\ar\ar \mbb{E}\left[\left|\int_0^t\int_{\mbb{R}} \sqrt{\gamma(\hat{\mu}_s, z)\hat{\mu}_s(z)}f(z) \hat{W}^{m_k}(d s, d z)\right|^2\right]\\
\ar\ar \quad \le K\mbb{E}\left[\int_0^t\int_{\mbb{R}} \hat{\mu}_s(z)J(z) d s d z\right] \le K\left\{\mbb{E}\left[\int_0^t\int_{\mbb{R}} \hat{\mu}_s(z)^2J(z) d s d z\right]\right\}^{1/2} < \infty.
\eeqnn 
It thus follows from \cite[Lemma~2.4]{XY19} that
\beqnn
 \<\hat{\mu}_t, f\> \ar=\ar \<\mu_0, f\> + \frac12\int_0^t \<\hat{\mu}_s, f''\>d s + \int_0^t \int_{\mbb{R}} \sqrt{\gamma(\hat{\mu}_s, z)\hat{\mu}_s(z)}f(z) \hat{W}(d s, d z).
\eeqnn
That completes the existence of solution to \eqref{1.1}. By Proposition~\ref{t2.4} one can see that $X_t(d x) = \mu_t(x)d x$ satisfies the MP (\ref{0.1}, \ref{0.2}), which implies the conclusion.
\qed

{\bf Proof of Theorem~\ref{t2}.} 
Assume that $(\mu_t(x))_{t \ge 0, x \in \mbb{R}}$ is a solution to SPDE \eqref{1.1}. Then it also satisfies \eqref{mu} by Proposition~\ref{t2.4}. For any $0 < r \le t$ and $p \ge 1$, we have
\beqnn
\mbb{E}\big[|\mu_t(x) - \mu_r(x)|^{2p}\big] \ar\le\ar K\Big|\int_\mbb{R} [p_t(x - z)- p_r(x - z)]\mu_0(z) d z\Big|^{2p}\\
\ar\ar + K\mbb{E}\left[\Big|\int_r^t\int_{\mbb{R}}p_{t - s}(x - z)\sqrt{\mu_s(z)}W(d s, d z)\Big|^{2p}\right]\\
\ar\ar+ K\mbb{E}\left[\Big|\int_0^r \int_{\mbb{R}}[p_{r - s}(x - z) - p_{t - s}(x - z)]\sqrt{\mu_s(z)}W(d s, d z)\Big|^{2p}\right]\\
\ar=: \ar K(I_1 + I_2 + I_3).
\eeqnn
Recall that $\mu_0$ is H{\"o}lder continuous with exponent $\lambda < \frac{1}{2}$. 
Similar with the proof of \cite[Theorem~1.1]{HLN13}, by H{\"o}lder's inequality  one can see that
\beqnn
I_1 \ar=\ar \Big|\int_\mbb{R} [p_t(x - z)- p_r(x - z)]\mu_0(z)d z\Big|^{2p}\\
\ar=\ar \Big|\mbb{E}\big[\mu_0(x + B_t) - \mu_0(x + B_r)\big]\Big|^{2p} \le K \mbb{E}\left[|B_{t - r}|^{2p\lambda}\right]\\
\ar\le\ar K\int_{\mbb{R}}|x|^{2p\lambda}\frac{1}{\sqrt{2\pi (t - r)}}e^{-\frac{x^2}{2(t - r)}} d x \le K(t - r)^{\lambda p},
\eeqnn
where $(B_t)_{t \ge 0}$ is a standard Brownian motion.  Moreover,  similar with \eqref{I_1} and \eqref{I_2} one can see that
\beqnn
I_2 \le K \mbb{E}\left[\Big|\int_r^t\int_{\mbb{R}}p_{t - s}(x - z)^2 \mu_s(z)d s d z\Big|^p\right]
\le K(t - r)^{1/4 + (p - 1)/2},
\eeqnn
and
\beqnn
I_3 \ar\le\ar K \mbb{E}\left[\Big|\int_0^r\int_{\mbb{R}}[p_{r - s}(x - z) - p_{t - s}(x - z)]^2\mu_s(z) d s d z\Big|^p\right]
\le K (t - r)^{\delta + (p -1)/2}
\eeqnn
with $\delta \in (0, 1/4).$ Then there exists  $\alpha \in (0, 1)$ such that
\beqlb\label{tr}
\mbb{E}\big[|\mu_t(x) - \mu_r(x)|^{2p}\big] \ar\le\ar K(t - r)^{\lambda p} + K(t - r)^{1/4 + (p - 1)/2} + K(t - r)^{\delta + (p -1)/2}\ccr
\ar\le\ar K(t - r)^{\alpha p/2}.
\eeqlb
Similarly, for any $x, y \in \mbb{R}$  we have
\beqlb\label{xy}
\mbb{E}\left[|\mu_t(x) - \mu_t(y)|^{2p}\right] \le K |x - y|^{2p\lambda} +  K|x - y|^{\beta p}
\eeqlb
with $\beta \in (0, 1).$ For $t, r \in [0, T]$ and $x, y \in \mbb{R},$ notice that
\beqnn
\mbb{E}\left[|\mu_t(x) - \mu_r(y)|^{2p}\right] \ar\le\ar K\mbb{E}\left[|\mu_t(x) - \mu_t(x)|^{2p}\right] + K\mbb{E}\left[|\mu_t(x) - \mu_t(y)|^{2p}\right].
\eeqnn
Combining with the above, \eqref{tr} and \eqref{xy}, the result follows from Kolmogorov's continuity criteria (see e.g. \cite[Corollary~1.2]{W86}).
\qed

\section{Proof of Theorem~\ref{main}}

In this section we show the weak uniqueness of the solution to the MP (\ref{0.1}, \ref{0.2}) under Condition \ref{0.4b}. The main idea of uniqueness is to relate the MP~(\ref{0.1}, \ref{0.2}) with a system of SPDEs, which is satisfied by a sequence of corresponding function-valued processes $(u_t^i)_{t \ge 0},\ i = 0, 1, \cdots, n.$ The weak uniqueness of the solution to the MP (\ref{0.1}, \ref{0.2}) follows from the pathwise uniqueness of the solution to the system of SPDEs.

For a solution $(X_t)_{t \ge 0}$ to the MP (\ref{0.1}, \ref{0.2}), we define the function-valued processes as
\beqlb\label{u}
 u_t^i(x) = X_t((a_i, x]), \qquad  a_i \le x < a_{i + 1}, \ i = 0, \cdots, n. 
\eeqlb
In Proposition~\ref{t1.2} we show that \eqref{u} is a solution to the following system of SPDEs:
\begin{numcases}{}
\begin{split}\label{2.1}
u_t^i(x) &= u_0^i(x) + \int_0^t \frac{1}{2}\Delta u_s^i(x) d s + \int_0^t\int_0^{u_s^i(x)}g_i(\nabla u^i_s(a_{i + 1}))W_i(d s, d z),\\
&\qquad\qquad \qquad\qquad \qquad \qquad\qquad x \in [a_i, a_{i + 1}),\ i = 0, \cdots, n - 1;
\end{split}
\\
\begin{split}\label{2.1a}
u_t^n(x) &= u_0^n(x) + \int_0^t \frac{1}{2}\Delta u_s^n(x) d s + \int_0^t\int_0^{u_s^n(x)}g_n(u_s^n(\infty))W_n(d s, d z),\\
& \qquad\qquad \qquad\qquad \qquad \qquad\qquad \qquad\qquad\qquad\qquad x \in [a_n, \infty).
\end{split}
\\
\begin{split}\label{0.7}
&u^i(a_i)=0,\;\; \;i=0,1,\cdots, n,\\
&\nabla u^{i-1}_t(a_i)=\nabla u^i_t(a_i),\;\;\;i=1,\cdots, n.
\end{split}
\end{numcases}
 where $u_s^n(\infty):=\lim_{x \rightarrow \infty}u_s^n(x)$ and $W_i(\d s, \d z),\ i = 0, \cdots, n$ are independent time-space white noises on $\mbb{R}_+ \times \mbb{R}_+$ with intensity $\d s \d z$. The pathwise uniqueness of the solution to the above system of SPDEs is obtained in Proposition \ref{t1.3}.

 The system of SPDEs (\ref{2.1}, \ref{2.1a}) can be understood in the following form:
 for any $\phi_i \in C_b^2[a_i, a_{i + 1}]$ with $\phi_i(a_i) = \phi'_i(a_{i + 1}) = 0,\ i = 0, 1, \cdots, n - 1,$
and $\phi_n \in C_b^2[a_n, \infty)$ with $\phi_n(a_n) = \phi_n(\infty) = 0$ (given $\phi_n(\infty) := \lim_{x \rightarrow \infty}\phi_n(x)$), we have
\begin{numcases}{}
\begin{split}\label{u^i}
\langle u_t^i, \phi_i\rangle &= \langle u_0^i, \phi_i\rangle + \frac{1}{2}\int_0^t \big[\langle u_s^i, \phi_i''\rangle + \phi_i(a_{i + 1})\nabla u^i_s(a_{i + 1})\big]d s\\
 & \qquad + \int_0^t\int_0^\infty\int_{a_i}^{a_{i + 1}} 1_{\{z \le u^i_s(x)\}}\phi_i(x)d xg_i(\nabla u^i_s(a_{i + 1}))W_i(d s, d z);
\end{split}
\\
\begin{split}\label{un}
\langle u_t^n, \phi_n\rangle &= \langle u_0^n, \phi_n\rangle + \frac{1}{2}\int_0^t \langle u_s^n, \phi_n''\rangle d s\\
 &\qquad + \int_0^t\int_0^\infty\int_{a_n}^\infty 1_{ \{z \le u^n_s(x)\}}\phi_n(x) d xg_n\big(u_s^n(\infty)\big)W_n(d s, d z).
\end{split}
\end{numcases}

\begin{prop}\label{t1.2}
Suppose that $(X_t)_{t \ge 0}$ is a solution to the MP (\ref{0.1}, \ref{0.2}). Then $\{u_t^i: t \ge 0, x \in [a_i, a_{i + 1})\}, i = 0, 1, \cdots, n$ defined as \eqref{u}, solves the group of SPDEs (\ref{2.1},\ref{2.1a}) with boundary condition \eqref{0.7}.
\end{prop}
\proof
 For any $\phi_i\in C^3_c(a_i,a_{i+1})$ with $i = 0, 1, \cdots, n-1$, by integration by parts, we have
 \beqnn
\<u^i_t,\phi_i'\> \ar=\ar - \<X_t,\phi_i\>=-M_t(\phi_i) - \<X_0,\phi_i\> - \frac12\int^t_0\<X_s,\phi_i''\>ds\\
\ar=\ar - M_t(\phi_i) + \<u^i_0,\phi_i'\> + \frac12\int^t_0\<u^i_s, \phi_i'''\> ds.
\eeqnn
Thus
\begin{equation}\label{eq0325b}
-M_t(\phi_i)=\<u^i_t,\phi_i'\> -\<u^i_0,\phi_i'\> -\frac12\int^t_0\<u^i_s, (\phi_i')''\> ds
\end{equation} 
is a continuous martingale. By Lemma~\ref{t2.3} we have
\beqnn
- M_t(\phi_i) = \int_0^t\int_{\mbb{R}}\phi_i(x)M(d s, d x)
\eeqnn
 with
\begin{eqnarray*}
\<-M(\phi_i)\>_t\ar=\ar \int^t_0g_i(\nabla u^i_s(a_{i + 1}))^2ds\int_{\mbb{R}} \phi_i(x)^2 X_s(d x)\\
\ar=\ar \int^t_0g_i(\nabla u^i_s(a_{i + 1}))^2d s \int_{0}^{u_s^i(a_{i + 1})}\phi_i(u_s^i(y)^{-1})^2 d y\\
\ar=\ar \int^t_0\int^\infty_0\(\int_{a_i}^{a_{i + 1}} 1_{\{y\le u^i_s(x)\}}\phi_i'(x)dx\)^2g_i(\nabla u^i_s(a_{i + 1}))^2ds dy,
\end{eqnarray*}
where $u_s^i(y)^{-1}$ denotes the generalized inverse of the nondecreasing function $u_s^i,$ that is,
\beqnn
u_s^i(y)^{-1} = \sup\{x \in [a_i, a_{i + 1}): u_s^i(x) \le y\}.
\eeqnn
Moreover, for $\phi_n \in C_c^3(a_n, \infty)$, one can see that
\begin{equation}\label{eq0325b}
-M_t(\phi_n) = \<u^n_t,\phi_n'\> -\<u^n_0,\phi_n'\> -\frac12\int^t_0\<u^n_s, \phi_n'''\> ds
\end{equation}
is a continuous martingale with
\beqnn
\<-M(\phi_n)\>_t\ar=\ar \int^t_0g_n(u_s^n(\infty))^2ds\int_{\mbb{R}} \phi_n(x)^2 X_s(d x)\\
\ar=\ar \int^t_0\int^{\infty}_0\(\int_{a_n}^{\infty} 1_{\{y \le u^n_s(x)\}}\phi_n'(x)dx\)^2g_n(u_s^n(\infty))^2ds dy.
\eeqnn
Similar with Lemma~\ref{t2.3}, the family $\{-M_t(\phi_i): t \ge 0, \phi_i \in C_c^3(a_i, a_{i + 1})\}$ determines a martingale measure $\{M_t(B): t \ge 0, B \in \mcr{B}(a_i, a_{i + 1})\}.$ Moreover, for $\phi_i \in C_c^3(a_i, a_{i + 1})$ and $\phi_j \in C_c^3(a_j, a_{j + 1})$ with $i \neq j,$ we have
\beqnn
\<-M(\phi_i), -M(\phi_j)\>_t \ar=\ar \int_0^t \d s \int_{\mbb{R}}\gamma(X_s, z)X_s(d z)\int_{\mbb{R}}\phi_i(x)\delta_z(d x)  \int_{\mbb{R}}\phi_j(y)\delta_z(d y)\\
\ar=\ar \int_0^t \d s \int_{\mbb{R}}\gamma(X_s, z)\phi_i(z)\phi_j(z)X_s(d z) = 0.
\eeqnn
By El Karoui and M\'{e}l\'{e}ard \cite[Theorem~III-7, Corollary~III-8]{EM90}, on some extension of the probability space on can define a sequence independent Gaussian white noise $W_i(d s, d u), i = 0, \cdots, n$ on $(0, \infty)^2$ based on $d s d u$ such that
\begin{eqnarray}\label{eq0325c}
-M_t(\phi_i)\ar=\ar\int^t_0\int^\infty_0\int_{a_i}^{a_{i + 1}} 1_{\{z \le u^i_s(x)\}}\phi_i'(x) \sqrt{\gamma(X_s, x)}dxW_i(d s, d z)\nonumber\\
\ar=\ar\int_{a_i}^{a_{i + 1}}\phi_i'(x) dx \int^t_0\int^{u^i_s(x)}_0g_i(\nabla u^i_s(a_{i + 1}))W_i(d s, d z)
\end{eqnarray}
for any $\phi_i \in C_c^3(a_i, a_{i + 1}), i = 0, \cdots, n - 1,$ and
\beqnn
-M_t(\phi_n)
\ar=\ar\int_{a_n}^{\infty} \phi_n'(x) dx\int^t_0\int^{u^n_s(x)}_0g_n( u^n_s(\infty))W_n(d s, d z)
\eeqnn
for any $\phi_n \in C_c^3(a_n, \infty).$ Then \eqref{u^i} and \eqref{un} follow from the approximation method and Lemma~\ref{l3.6} in Appendix.  
\qed

\begin{lemma}\label{l3.2}
Suppose that $(u_t^n)_{t \ge 0}$ satisfies \eqref{2.1a}.
Then $(u_t^n(\infty))_{t \ge 0}$ satisfies
 \beqlb\label{2.10}
u_t^n(\infty)=u_0^n(\infty)+\int_0^t
\int_0^{u_s^n(\infty)}g_n(u_s^n(\infty))W_n(d s, d z).
 \eeqlb
\end{lemma}

\proof
Recall that $p_t(x)=\frac{1}{\sqrt{2\pi}t}\e^{-x^2/(2t)}$ and let
\beqnn
q_t^x(y):=p_t(x+a_n-y)-p_t(x-a_n+y)
 \eeqnn
for $t>0$ and $x,y\ge a_n$.
Then \eqref{2.1a} can be written into the following mild form:
 \beqlb\label{2.8}
u_t^n(x)\ar=\ar\<u_0^n,q_t^x\> \nonumber\\
\ar\ar +\int_0^t\int_0^\infty
\Big[\int_{a_n}^\infty1_{\{z \le u_s^n(y)\}}q_{t-s}^x(y)\d y\Big]g_n(u_s^n(\infty)) W_n(d s, d z) 
 \eeqlb
for any $x \ge a_n.$ By a change of variable, we have
 \beqnn
\<u_0^n,q_t^x\>
 \ar=\ar
\int_{a_n}^\infty u_0^n(y)[p_t(x+a_n-y)-p_t(x-a_n+y)]d y \\
 \ar=\ar
\int_{-\infty}^x u_0^n(x+a_n-z)p_t(z)d z
-\int_x^\infty u_0^n(z-x+a_n)p_t(z)d z \\
 \ar\to\ar
u_0^n(\infty)\int_{-\infty}^\infty p_t(z)d z=u_0^n(\infty)
 \eeqnn
and
 \beqnn
\int_{a_n}^\infty1_{\{z \le u_s^n(y)\}}q_{t-s}^x(y) d y\to1_{\{z \le u_s^n(\infty)\}}
 \eeqnn
as $x\to\infty,$ which ends the proof.
\qed

\begin{prop}\label{t1.3}(Pathwise uniqueness)
 Suppose that $(u_t)_{t \ge 0}$ and $(\tilde{u}_t)_{t \ge 0}$ are two solutions to (\ref{2.1}, \ref{2.1a}) satisfying the boundary conditions~\eqref{0.7}. If $u_0(x) = \tilde{u}_0(x)$ for all $x \in \mbb{R},$ then $\mbb{P}\{u_t(x) = \tilde{u}_t(x)\ \text{for all}\ t \ge 0\ \text{and}\ x \in \mbb{R}\} = 1.$ 
\end{prop}

\proof
The pathwise uniqueness of the solution to \eqref{2.10} holds by \cite[Theorem 2.1]{DL12}. Moreover, the pathwise uniqueness of the solution $\{u_t^n(x): t \ge 0, x \in [a_n, \infty)\}$ holds for \eqref{2.1a} by \cite[Theorem~1.4]{XY20}. That implies the strong uniqueness of $(\nabla u_t^n(a_n))_{t \ge 0}.$ By the proof of \cite[Theorem~1.4]{XY20} and induction method, one can get the pathwise uniqueness of the solution to \eqref{2.1} for $i = 0, 1, \cdots, n-1.$ The proof ends here.
 \qed

{\bf Proof of Theorem~\ref{main}}
The result is a direct conclusion of Proposition~\ref{t1.3}.
\qed

\section{Appendix}


In this section we give some results about the process $(u_t(x))_{t \ge 0, x \in [0, 1]}$ satisfying the following SPDE:
\beqnn
u_t(x) = u_0(x) + \int_0^t \frac{1}{2}\Delta u_s(x) d s + \int_0^t\int_0^{u_s(x)}g(\nabla u_s(1))W(d s, d z),
\eeqnn
where $g$ is a positive continuous bounded function from $\mbb{R}_+$ to $\mbb{R}_+$, and $W(d s, d z)$ is a time-space Gaussian white noise with density $d s d z$. Let $\Phi\in C_c^2(0,1)$ satisfies $0\le \Phi\le 2$ and $\int_{0}^1\Phi(x)\d x=1$. 
For $k\ge1$ and $x\in[0,1]$ let
 \beqlb\label{hn}
h_k(x):= \int_0^{kx}\Phi(z)\d z\cdot \int_{x^k}^1\Phi(z)\d z.
 \eeqlb
Then $h_k\in C_c^2(0,1)$ for all $k\ge1$.

\begin{lemma}\label{l3.5}
Suppose that $f\in C[0,1]$ with $f'(1)$ and $f'(0)$ exist.
Then
 \beqnn
\lim_{k\to\infty}\<f,h_k'\> =  f(0)-f(1),\quad
\lim_{k\to\infty}\<f,h_k''\> = f'(1)-f'(0).
 \eeqnn
\end{lemma}
\proof
Observe that for each $n\ge1$,
 \beqnn
h_k'(x)=k\Phi(kx)
\int_{x^k}^1\Phi(z) d z
-kx^{k-1}\Phi(x^k)\int_0^{kx}\Phi(z)d z,\qquad x\in[0,1].
 \eeqnn
Then by change of variables and dominated convergence,
as $k\to\infty$,
 \beqnn
\<f,h_k'\>
 \ar=\ar
\int_0^1f(x)k\Phi(kx)\Big[\int_{x^k}^1\Phi(z) d z\Big] d x
-\int_0^1f(x)kx^{k-1}\Phi(x^k)\Big[\int_0^{kx}\Phi(z) d z\Big] d x \\
 \ar=\ar
\int_0^1f(y/k)\Phi(y)\Big[\int_{(y/k)^k}^1\Phi(z) d z\Big] d y
-\int_0^1f(y^{1/k})\Phi(y)\Big[\int_0^{ky^{1/k}}\Phi(z) d z\Big] d y
 \eeqnn
converges to $f(0)-f(1)$,
which gives the first assertion.

In the following we prove the second assertion.
Observe that
 \beqlb\label{5.1}
h_k''(x)
 \ar=\ar
k^2\Phi'(kx)\int_{x^k}^1\Phi(z)\d z-2k^2 x^{k-1}\Phi(kx)\Phi(x^k) \nonumber\\
\ar\ar 
-\big[k(k-1)x^{k-2}\Phi(x^k)
+k^2x^{2k-2}\Phi'(x^k)\big]\int_0^{kx}\Phi(z) d z \nonumber\\
 \ar=:\ar
M_{1,k}(x)-2M_{2,k}(x) - M_{3,k}(x).
 \eeqlb
By change of variables and dominated convergence again,
as $k\to\infty$,
 \beqlb\label{5.2}
\int_0^1[f(x)-f(0)]M_{1,k}(x) d x
 \ar=\ar
\int_0^1k[f(y/k)-f(0)]\Phi'(y)\Big[\int_{y^kk^{-k}}^1\Phi(z) d z\Big] d y  \nonumber\\
 \ar\to\ar
f'(0)\int_0^1y \Phi'(y) d y
=-f'(0)
 \eeqlb
and
 \beqlb\label{5.3}
\int_0^1[f(x)-f(1)]M_{2,k}(x) d x
 \ar=\ar\int_0^1\frac{f(y^{1/k})-f(1)}{y^{1/k}-1}k(y^{1/k}-1)\Phi(y^{1/k})\Phi(y) d y\nonumber\\
\ar\to\ar 0.
 \eeqlb
Similarly, as $k\to\infty$,
 \beqlb\label{5.4}
 &\ &
\int_0^1[f(x)-f(1)]M_{3,k}(x) d x \nonumber\\
 &\ &\quad=
\int_0^1\frac{f(y^{1/k})-f(1)}{y^{1/k}-1}k(y^{1/k}-1)
\big[k^{-1}(k-1)y^{-1/k}\Phi(y) \nonumber\\
 &\ &\quad\qquad
+y^{(k-1)/k}\Phi'(y)\big]\cdot\Big[\int_0^{ky^{1/k}}\Phi(z) d z\Big] d y\nonumber\\
 &\ &\quad\to
f'(1)\int_0^1 \ln y[\Phi(y)+y\Phi'(y)] d y=-f'(1).
 \eeqlb
Applying integration by parts and the fact $0\le\Phi\le 2$ and $\mbox{supp}(\Phi)\subset(0,1)$,
 \beqnn
\int_0^1M_{1,k}(x) d x
 \ar=\ar
\int_0^1\Big(\int_0^{kx}\Phi(z) d z\Big)''\cdot \Big(\int_{x^k}^1\Phi(z) d z\Big) d x=\int_0^1M_{2,k}(x) d x \\
 \ar=\ar
k\int_0^1 \Phi(ky^{1/k})\Phi(y) d y
=
k\int_0^{k^{-k}} \Phi(ky^{1/k})\Phi(y) d y
\le 4 k^{1-k}
 \eeqnn
and
 \beqnn
\int_0^1M_{3,k}(x) d x \ar=\ar
-\int_0^1\Big(\int_{x^k}^1\Phi(z) d z\Big)''\cdot\Big(\int_0^{kx}\Phi(z) d z\Big) d x\\
\ar=\ar \int_0^1M_{2,k}(x) d x
\le 4 k^{1-k}.
 \eeqnn
Then combining \eqref{5.1} with \eqref{5.2}-\eqref{5.4}
one completes the proof.
\qed

\begin{lemma}\label{l3.6}
Suppose that for each $\phi \in C_c^2(0, 1)$,  $(u_t)_{t\ge0}$ satisfies
 \beqlb\label{uiphi}
\<u_t,\phi\>
\ar=\ar
\<u_0,\phi\>  + \frac12\int_0^t\<u_s,\phi''\> d s\nonumber\\
&\ & +\int_0^t\int_0^\infty g (\nabla u_s(1)) \Big[\int_0^1 1_{\{z \le u_s (x)\}}\phi(x)d x\Big]W (d s, d z).
 \eeqlb
Then for each $\phi\in C_b^2[0, 1]$,
 \beqlb\label{3.5}
\<u_t,\phi\>
\ar=\ar
\<u_0,\phi\>  + \frac12\int_0^t\big[\<u_s,\phi''\>
+F_s(\phi)] d s\nonumber\\
&\ & + \int_0^t\int_0^\infty g(\nabla u_s(1))\Big[\int_0^11_{\{z \le u_s(x)\}}\phi(x)\d x\Big]W(d s, d z),
 \eeqlb
where
 \beqlb\label{3.9}
F_s(\phi):=[\phi(1)\nabla u_s(1)-\phi(0)\nabla u_s(0)]-[u_s(1)\phi'(1)- u_s(0)\phi'(0)].
 \eeqlb
\end{lemma}

\proof
Recall $h_k$ in \eqref{hn}.
For $m\ge1$ define stopping time $\tau_m$ by
 \beqnn
\tau_m:=\inf\Big\{t\ge0:\sup_{x\in[0,1]}|u_t(x)|\ge m\Big\}
 \eeqnn
with the convention $\inf\emptyset=\infty$.
Then $\lim_{m \to\infty}\tau_m=\infty$ almost surely.
It follows from \eqref{uiphi} that
\beqlb\label{3.3}
\<u_{t\wedge\tau_m},\phi h_{k}\>
 \ar=\ar \int_0^{t\wedge\tau_m}\int_0^\infty g(\nabla u_s(1))\Big[\int_0^11_{\{z \le u_s(x)\}}\phi(x)h_{k}(x)d x\Big]W(d s, d z)\nonumber\\
 &\ &
+\<u_0,\phi h_{k}\> +\frac12\int_0^{t\wedge\tau_m}\<u_s,(\phi h_{k})''\> d s .
 \eeqlb
Notice that
 \beqnn
\<u_s,(\phi h_{k})''\>
=
\<u_s,\phi''h_{k}\>
+2\<u_s,\phi'h'_{k}\>
+\<u_s,\phi h''_{k}\>.
 \eeqnn
It follows from Lemma \ref{l3.5} that
 \beqnn
\lim_{k\to\infty}\<u_s,(\phi h_{k})''\>
 \ar=\ar
\<u_s,\phi''\>
+[u_s(0)\phi'(0)-u_s(1)\phi'(1)]
-[\phi(0)\nabla u_s(0)-\phi(1)\nabla u_s(1)] \\
 \ar=\ar
\<u_s,\phi''\> + F_s(\phi).
 \eeqnn
Thus letting $k\to\infty$ in \eqref{3.3} we obtain
 \beqnn
\<u_{t\wedge\tau_m},\phi\>
\ar=\ar
\<u_0,\phi\>
+\frac12\int_0^{t\wedge\tau_m}\big[\<u_s,\phi''\>
+ F_s(\phi)\big]d s\\
&\ & +\int_0^{t\wedge\tau_m}\int_0^\infty g(\nabla u_s(1))\Big[\int_0^1 1_{\{z \le u_s(x)\}}\phi(x)d x\Big]W(d s, d z).
 \eeqnn
Letting $m\to\infty$ the result holds.
\qed

\medskip

{\bf Acknowledgements} The research of L. Ji was supported in part by the fellowship of China Postdoctoral Science Foundation  2020M68194; the research of J. Xiong  was supported in part by NSFC grants 61873325, 11831010 and SUSTech fund Y01286120; the research of X. Yang was supported in part by NSFC grants 11771018 and 12061004.


\begin{thebibliography}{999}

\bibitem{D75} Dawson, D. A. (1975): Stochastic evolution equations and related measure processes. \textit{J. Multivariate Anal.}, {\bf 5}(1):1--52.

\bibitem{DL03} Dawson, D. A. and Li, Z. (2003): Construction of immigration superprocesses with dependent spatial motion from one-dimensional excursions, \textit{Probab. Theory Related Fields} 127, 37--61.

\bibitem{DL12}  Dawson, D. A. and Li, Z. (2012): Stochastic equations, flows and measure-valued processes. \textit{Ann. Probab.} 40(2): 813--857.


\bibitem{DK99} Donnelley, P. and Kurtz, T. (1999): Particle representations for measure-valued population models, \textit{Ann. Probab.} 27(1): 166--205.

\bibitem{EM90} N. El Karoui, S. M\`{e}l\`{e}ard (1990): Martingale measures and stochastic calculus. \textit{Probab. Theory Related Fields} 84, 83--101.


\bibitem{EK86} Ethier, S. N. and Kurtz, T. G. (1986): \textit{Markov Processes: Characterization and Convergence.} Wiley, New York.

\bibitem{FL04}  Fu, Z. and Li, Z. (2004): Measure-valued diffusions and stochastic equations with Poisson process. \textit{Osaka J. Math.} {\bf 41}(3): 727--744.

\bibitem{HLN13} Hu, Y., Lu, F. and Nualart, D. (2013): H\"{o}lder continuity of the solutions for a class of nonlinear SPDE's arising from one dimensional superprocesses. \textit{Probab. Thery Relat. Fields} {\bf 156}: 27--49.

\bibitem{HLY14} He, H., Li, Z. and Yang, X. (2014): Stochastic equations of super-L\'{e}vy processes with general branching mechanism. \textit{Stochastic Process. Appl.} {\bf 124}(4), 1519--1565.



\bibitem{K02} Kallenberg, O. (2002): \textit{Foundations of modern probability}. Second edition. Probability and its Applications (New York). Springer-Verlag, New York.




\bibitem{L11} Li, Z. (2011): \textit{Measure-valued branching Markov processes.} Springer.


\bibitem{LWX05} Li, Z., Wang, H. and Xiong, J. (2005): Conditional log-Laplace functionals of immigration superprocesses with dependent spatial motion. \textit{Acta Appl. Math.} {\bf 88}, 143-175.


\bibitem{M85} Mitoma, I. (1985): An $\infty$-dimensional inhomogeneous Langevin's equation. {\em J. Funct. Anal.} {\bf 61}: 342--359.

\bibitem{MX15} Mytnik, L. and Xiong, J. (2015): Well-posedness of the martingale problem for superprocess with interaction. \textit{Illinois J. Math.} {\bf 59} (2): 485--497.

\bibitem{P95} Perkins E. (1995): \textit{On the martingale problem for interactive measure-valued branching diffusions.} Mem. Amer. Math. Soc. 115.

\bibitem{P02} Perkins, E. (2002): \textit{Dawson-Watanabe superprocesses and measure-valued diffusions}. Lectures on probability theory and statistics (Saint-Flour, 1999), Lecture Notes in Math., vol. 1781. Springer, Berlin.


\bibitem{S90} Shiga, T. (1990): A stochastic equation based on a Poisson system for a class of measure valued diffusion processes \textit{J. Math. Kyoto Univ.} {\bf 30}, 245--279.

\bibitem{WillettWong} Willett, Z. and Wong, J. (1965): On the discrete analogues of some generalizations of Gronwall's inequality. \textit{Monatsh. Math.} {\bf 69}, 362--367.


\bibitem{W86} Walsh, J. (1986): An introduction to stochastic partial differential equations. \textit{Lecture Notes in Math.} {\bf 1180}, 266--439. Spring, Berlin.


\bibitem{X13a} Xiong, J. (2013): \textit{Three classes of nonlinear stochastic partial differential equations}. World Scientific, Hackensack, NJ.

\bibitem{X13b} Xiong, J. (2013): Super-Brownian motion as the unique strong solution to an SPDE. \textit{Ann. Probab.} {\bf 41}(2): 1030--1054.

\bibitem{XY16} Xiong, J. and Yang, X. (2016): Superprocesses with interaction and immigration. \textit{Stochastic Process. Appl.} {\bf 126}(11): 3377--3401.


\bibitem{XY19} Xiong, J. and Yang, X. (2019): Existence and pathwise uniqueness to an SPDE driven by $\alpha$-stable colored noise. \textit{Stochastic Process. Appl.} {\bf 129}: 2681--2722.

\bibitem{XY20} Xiong, J. and Yang, X. (2021+): SPDEs with non-Lipschitz coefficients and nonhomogeneous boundary conditions. arXiv: 2006.01009.

\end{thebibliography}
\end{document}